\documentclass[10pt,a4paper]{article}
\usepackage{amsfonts}
\usepackage{}
\usepackage{amsmath}
\usepackage{latexsym}
\usepackage{amssymb}
\usepackage{amstext}
\usepackage{fancyhdr}
\usepackage{indentfirst}
\usepackage{graphicx}
\usepackage{epsfig}
\usepackage{pgf}
\usepackage{listings}
\usepackage{cases}
\usepackage{epstopdf}
\usepackage{graphicx}
\usepackage{enumerate}
\usepackage{amsthm}
\usepackage{mathrsfs}
\usepackage{abstract}
\usepackage{caption}
\date{}

\begin{document}
\title{The signless Laplacian state transfer in $\mathcal{Q}$-graph}
\author{Xiao-Qin Zhang$^a$, Shu-Yu Cui$^b$\footnote{Corresponding author. E-mail address: cuishuyu@zjnu.cn or cuishuyu@163.com.}, Gui-Xian Tian$^{a*}$\\%EndAName
    {\small{\it $^a$Department of Mathematics,}}
    {\small{\it Zhejiang Normal University, Jinhua, Zhejiang, 321004, P.R. China}}\\
    {\small{\it $^b$Xingzhi College, Zhejiang Normal University, Jinhua, Zhejiang, 321004, P.R. China}}
}\maketitle
\begin{abstract}
The $\mathcal{Q}$-graph of a graph $G$, denoted by $\mathcal{Q}(G)$, is the graph derived from $G$ by plugging a new vertex to each edge of $G$ and adding a new edge between two new vertices which lie on adjacent edges of $G$. In this paper, we consider to study the existence of the signless Laplacian perfect state transfer and signless Laplacian pretty good state transfer in $\mathcal{Q}$-graphs of graphs. We show that, if all the signless Laplacian eigenvalues of a regular graph $G$ are integers, then the $\mathcal{Q}$-graph of $G$ has no signless Laplacian perfect state transfer. We also give a sufficient condition that the $\mathcal{Q}$-graph of a regular graph has signless Laplacian pretty good state transfer when $G$ has signless Laplacian perfect state transfer between two specific vertices.


\end{abstract}

\section{Introduction}

In 1998, the concept of continuous-time quantum walk on graphs was proposed firstly by Farhi and Gutman in \cite{Farhi}. Suppose that $M$ is a symmetric matrix associated with a connected graph $G$, such as adjacency matrix, Laplacian matrix, signless Laplacian matrix of $G$ and so on. Then the \emph{transition matrix} of continuous-time quantum walk on $G$ is defined by the unitary matrix
\begin{equation}
H_{M}(t)=\text{exp}(-itM),
\end{equation}
where $t\in \mathbb{R}$ and $i^2=-1$.
Bose \cite{Bose} in $2003$ studied the problem of information transition in a quantum spin system. Subsequently,
Christandl et al. \cite{Christandl} proved that this problem can be reduced to the perfect state transfer in quantum walk. Let $e_{u}^n$ denote the characteristic vector of order $n$ corresponding to the vertex $u$ of $G$. In the absence of confusion, $e_{u}^n$ is abbreviated as $e_{u}$. For two vertices $u$ and $v$ of $G$, if
\begin{equation}
\text{exp}(-itM)e_{u}=\lambda e_{v}
\end{equation}
with $ |\lambda| =1$, then $G$ is said to have \textit{perfect state transfer} relative to matrix $M$ between vertices $u$ and $v$ at time $t$. This concept has played a crucial role in quantum information and quantum algorithm. Characterizing graphs having perfect state transfer has been enjoyed world-wide attention among physic and mathematics communities. However, it is well known that graphs having perfect state transfer are rare. Godsil in \cite{Godsil2012} proposed a new concept called pretty good state transfer, whose condition is more relaxing in comparison to perfect state transfer. A graph $G$ is said to have \emph{pretty good state transfer} between vertices $u$ and $v$ at some time $t$ if
\begin{equation}
|\text{exp}(-itM)_{uv}|>1-\epsilon,
\end{equation}
for any $\epsilon>0$.

Many papers have focused mostly on quantum walk relative to adjacency matrix and Laplacian matrix and obtained some excellent results. For example, in the adjacency matrix case, Godsil in \cite{equation} proved that, for any $k\in\mathbb{Z}^+$, there exist at most finitely graphs with maximum degree $k$ admitting perfect state transfer. Bose et al. in \cite{S.Bose} showed the complete graph $K_{n}$ has no perfect state transfer, but $K_{n}$ with a missing edge has perfect state transfer between the two vertices of this missing edge for positive integer $n$. Ge et al. in \cite{Ge Yang} described some new constructions of perfect state transfer graphs using variants of the double cones and graph products, such as weak and lexicographic products, irregular and glued double cones and so on. Coutinho and Godsil in \cite{Coutinho2016} also constructed many new graphs of perfect state transfer using graph products and double covers of graphs. Ackelsberg et al. in \cite{Ackelsberg2017} showed that the corona of two graphs has no perfect state transfer under suitable conditions, but it has pretty good state transfer in some special cases. For more details on this area, we refer the reader to \cite{Coutinho2014,Godsil2012} and the cited references therein.

In the Laplacian matrix case, Alvir et al in \cite{Alvir2016} studied the
perfect state transfer in Laplacian quantum walks. It was indicated that quantum walks based on the adjacency matrix, Laplacian matrix and signless Laplacian matrix are all equivalent for a regular graph. They also investigated the Laplacian perfect state transfer on graph joins and provided a characterization of Laplacian perfect state transfer on the double cones. Coutinho and Liu in \cite{Coutinho2015} proved that a tree of order $n\geq 3$ has no Laplacian perfect state transfer. However, Banchi et al. in \cite{Banchi2017} showed that a path of order $n$ has Laplacian pretty good state transfer if and only if $n$ is a power of 2. Recently, Ackelsberg et al. in \cite{Ackelsberg2016} showed that the corona of two graphs has no Laplacian perfect state transfer, but it has Laplacian pretty good state transfer with some mild conditions. Li and Liu in \cite{Yipeng Li} showed  $\mathcal{Q}$-graph of an $r$-regular graph $G$ has no Laplacian perfect state transfer when $r+1$ is a prime number, but it has Laplacian pretty good state transfer. Liu and Wang in \cite{Xiaogang Liu} proved that total graph of an $r$-regular graph $G$ has no Laplacian perfect state transfer when $r+1$ is not a Laplacian eigenvalue of $G$. They also gave a sufficient condition on Laplacian pretty good state transfer for total graph of a regular graph. For more information about Laplacian state transfer of graphs, readers may refer to \cite{Alvir2016,Coutinho2014,Wang2021} and the cited references therein.

Compared with the research on quantum state transfer relative to adjacency matrix and Laplacian matrix of graphs, there have been far less studies on the problem of signless Laplacian state transfer of graphs. For example, Alvir et al in \cite{Alvir2016} proved that, if $H$ is a $(\frac{n}{2}-1)$-regular graph of order $n\geq2$, then the double cones $\overline{K}_2+H$ has signless Laplacian perfect state transfer. It is proved \cite{Alvir2016} that there exists somewhat interesting tie between
quantum walks relative to the signless Laplacian matrix of a graph $G$ and adjacency matrix of its line graph $l(G)$. With the help of this, Alvir et al \cite{Alvir2016} proved that the path $P_n$ has no signless Laplacian prefect state transfer for $n\geq5$. Recently, Kempton et al. in \cite{potential,involution} mainly investigated perfect state transfer and pretty good state transfer with potential. Note that, if the potential on each vertex of graphs equals its degree, then it will result in the signless Laplacian quantum state transfer of graphs. In $2021$, Tian et al. \cite{Tian2021} showed that there is no signless Laplacian perfect state transfer on corona graph $G\circ K_{m}$ if $m$ equals one or a prime number. They also showed, if $m$ is an even number, then $K_{2}\circ \overline{K}_{m}$ has no signless Laplacian perfect state transfer between two vertices of $K_{2}$. However, $G\circ \overline{K}_{m}$ has signless Laplacian pretty good state transfer under some special conditions. For some properties about signless Laplacian matrices and applications in computer science, see \cite{Cvetkovic2010,Cvetkovic2009,Cvetkovic2011,Dam2003} and the cited references therein.

Motivated by above some results, we mainly consider to study on perfect state transfer and pretty good state transfer relative to signless Laplacian matrix in $\mathcal {Q}$-graph. The $\mathcal {Q}$-\textit{graph} \cite{Cvetkovic2010} of a graph $G$, denoted by $\mathcal {Q}(G)$, is the graph obtained from $G$ by inverting one new vertex in each edge and joining these new vertices by edges which lie on adjacent edge in $G$. In our work, we show that there is no signless Laplacian perfect state transfer on $\mathcal {Q}$-graph of a regular graph with a mild condition. Furthermore, we also present a sufficient condition for the $\mathcal {Q}$-graph admitting signless Laplacian pretty good state transfer.

\section{Preliminaries}

For a graph $G$ of order $n$, let $A(G)$ and $D(G)$ denote its adjacent matrix and degree diagonal matrix, respectively. Then $Q_{G}=A(G)+D(G)$ is the signless Laplacian matrix of $G$. The eigenvalues of $Q_{G}$ are called the \textit{signless Laplacian eigenvalues} of $G$. Let $q_{0}>q_{1}>\cdots>q_{d}$ be all distinct eigenvalues of $Q_{G}$ and $x_{1}^{(i)}, \ldots, x_{l_{i}}^{(i)}$ be the unit orthogonal eigenvectors corresponding to eigenvalue $q_{i}$ with multiplicity $l_{i}$, $i=0,1, \ldots, d$.  Let $x^T$ denote the transpose of a column vector $x$. Then, the eigenprojector corresponding to eigenvalue $q_{i}$ can be written as the following matrix:
\begin{equation}
{f_{{q_i}}} = \sum\limits_{j = 1}^{{l_i}} {x_j^{(i)}(x_j^{(i)}} {)^T},
\end{equation}
and $\sum_{i=0}^{d} f_{{q_i}} =I_n$, where $I_n$ denotes the identity matrix of order $n$. It is easy to find  that the $f_{{q_i}} $ is idempotent matrix $(f_{{q_i}}^{2}=f_{{q_i}})$ and $f_{{q_i}}f_{{q_k}}=0$ for $i\neq k$. We can obtain the following alternative expression of the signless Laplacian matrix $Q_{G}$ by the eigenprojector:
\begin{equation}
Q_{G}=Q_{G}\sum_{i=0}^{d} f_{{q_i}}=\sum_{i=1}^{d}\sum_{j=1}^{l_{i}} Q_{G} x_{j}^{(i)} (x_{j}^{(i)})^T=\sum_{i=0}^{d}\sum_{j=1}^{l_{i}} q_{i} x_{j}^{(i)} (x_{j}^{(i)})^T=\sum_{i=1}^{d}q_{i} f_{{q_i}},
\end{equation}
which is called the \textit{spectral decomposition} of $Q_{G}$. Since the continuous-time quantum walk on graph $G$ relative to the signless Laplacian matrix $Q_{G}$ is the unitary matrix
\begin{equation}
H_{Q_{G}}(t)=\text{exp}(-itQ_{G}),
\end{equation}
then
\begin{equation}\label{7}
H_{Q_{G}}(t)=\text{exp}(-itQ_G)=\sum_{k\geq 0}\frac{(-i)^{k}Q_{G}^{k}t^{k}}{k!}=\sum_{k\geq 0}\frac{(-i)^{k}\sum_{i=0}^{i=d}q_{i}^{k} f_{{q_i}}t^{k}}{k!}=\sum_{i=1}^{i=d} \text{exp}(-it q_{i})f_{{q_i}}.
\end{equation}

For a vertex $u$ of a graph $G$, denote its \textit{signless Laplaian eigenvalue support} in graph $G$ by $\text{supp}_{Q_{G}}(u)$, which is the set of all the eigenvalues of $Q_{G}$ satisfying $f_{{q_i}}e_{u}\neq 0$, where $e_{u}$ is called the \textit{characteristic vector} of $u$ (the $u$-th entry of the column vector $e_{u}$ is one, otherwise zero). The vertices $u$ and $v$ are said to be \textit{strongly signless Laplacian cospectral} if $f_{q}e_{u}=\pm f_{q}e_{v}$ for all eigenvalues $q$ of $Q_{G}$. Denote the set of all the eigenvalues such that $f_{q}e_{u}=f_{q}e_{v}$ by $S^{+}$, and denote the set of all the eigenvalues such that $f_{q}e_{u}=-f_{q}e_{v}$ by $S^{-}$.

The following are some main theorems and lemmas which will help us to study singless Laplacian perfect state transfer or singless Laplacian pretty good state transfer of $\mathcal {Q}$-graph.

\paragraph{Theorem 2.1.}(Coutinho \cite{Coutinho2014}) Let $G$ be a graph of order $n\geq2$ with the vertex set $V(G)$, and $u, v\in V(G)$.
Suppose that $q_{0}$ is the maxmum signless Laplacian eigenvalue in $G$.
Then $G$ admits signless Laplacian perfect state transfer between the vertices $u$ and  $v$ if and only if the following conditions hold.
\begin{enumerate}[(i)]
 \item Two vertices $u$ and $v$ are strongly cospectral relative to signless Laplacian.
 \item Non-zero elements in $\text{supp}_{Q_G}(u)$ are either all integers or all quadratic integers. Moreover, there exists a square-free integer $\Delta$ and integers $a$, $b_q$ such that, for each signless Laplacian eigenvalue $q\in \text{supp}_{Q_G}(u)$,
  \begin{equation*}
  q=\frac{1}{2}(a+b_q\sqrt{\Delta}).
  \end{equation*}
  Here we allow $\Delta=1$ if all signless Laplacian eigenvalues in $\text{supp}_{Q_G}(u)$ are integers, and $a=0$ if all signless Laplacian eigenvalues in $\text{supp}_{Q_G}(u)$ are all multiples of $\sqrt{\Delta}$.
 \item $q\in{S^+}$ if and only if $\dfrac{q_0-q}{g\sqrt{\Delta}}$ is even, where
 $$g=\gcd\left(\left\{\dfrac{q_0-q}{\sqrt{\Delta}}:q\in \text{supp}_{Q_G}(u)\right\}\right).$$
 %Then \subitem{(a)} and
     %\subitem{(b)}  $q\in{S^-}$ if and only if $\dfrac{q_0-q}{g\sqrt{\Delta}}$ is odd.
\end{enumerate}
Moreover, if the above conditions hold, then the following also hold.
\begin{enumerate}[(1)]
     \item There exists a minimum time $\tau_0>0$ at which signless Laplacian perfect state transfer occurs between $u$ and $v$, and
     \begin{equation*}
     \tau_0=\frac{1}{g}\dfrac{\pi}{\sqrt{\Delta}}.
     \end{equation*}
     \item The time of signless Laplacian perfect state transfer $\tau$ is an odd multiple of $\tau_0$.
     \item The phase of signless Laplacian perfect state transfer is given by $\lambda=e^{-itq_0}$.
\end{enumerate}

\paragraph{Theorem 2.2.}(Hardy and Wright \cite{numbertheory}) Assume that $1, q_1, \ldots, q_m$ are linearly independent over the set of rational number $\mathbb{Q}$. Then, for any real numbers $\alpha_1, \ldots, \alpha_m$ and positive real numbers $N$, $\epsilon$, there exist integers $l>N$ and $\gamma_1, \ldots, \gamma_m$ such that
\begin{equation}
\vert{lq_{k}-\gamma_k-\alpha_k}\vert<\epsilon,
\end{equation}
for each $k=1,\ldots,m$.

\paragraph{Lemma 2.3.}(Richards \cite{galoistheory}) The set $\{\sqrt{\Delta}:\Delta$ is a square-free integer$\}$ is linearly independent over
$\mathbb{Q}$.

\paragraph{Lemma 2.4.}(Coutinho \cite{Coutinho2015}) A real number
$\lambda$ is a quadratic integer if and only if there exist integers $a$, $b$ and $\Delta$ such that $\Delta$ is square-free and one of the following cases holds.
\begin{enumerate}[(i)]
    \item $\lambda=a+b\sqrt{\Delta}$ and $\Delta\equiv2,3\;(\text{mod}\;4)$.
    \item $\lambda=\frac{1}{2}(a+b\sqrt{\Delta}),\; \Delta\equiv1\;(\text{mod}\;4)$, and either $a$ and $b$ are both even or both odd.
\end{enumerate}

\section{Signless Laplacian eigenvalues and eigenprojectors of $\mathcal {Q}$-graph}

Let $G$ be a graph of order $n$ with vertex set $V(G)=\{v_{1}, v_{2}, \ldots, v_{n}\}$ and edge set $E(G)=\{e_{1},e_{2}, \ldots, e_{m}\}$. The \textit{vertex-edge incidence matrix} of $G$ is an $n\times m$ matrix $R_{G}=(r_{ij})_{n\times m}$, in which $r_{ij}=1$ if the vertex $v_{i}$ is incident to the edge $e_{j}$, otherwise $r_{ij}=0$. Let $y_{1},y_{2},\ldots,y_{l}$ be all unit orthogonal vectors such that $R_Gy=0$. It is well known \cite{Cvetkovic1995} that $l=m-n$ if $G$ is non-bipartite graph and $l=m-n+1$ if $G$ is bipartite graph. Denote the column vector with all of its entries one by $j_n$.

Next, we first give the signless Laplacian eigenvalues and corresponding signless Laplacian eigenvectors of $\mathcal {Q}(G)$.

\paragraph{Theorem 3.1.} Let $G$ be an $r$-regular non-bipartite connected graph of order $n$, with $m$ edges and $r\geq2$. Also let $2r=q_{0}>q_{1}>\cdots>q_{d}$ be all distinct signless Laplacian eigenvalues of $G$ and $x_{1}^{(i)},x_{2}^{(i)}, \ldots, x_{l_{i}}^{(i)}$ be the unit orthogonal eigenvectors corresponding to eigenvalue $q_{i}$ with multiplicity $l_{i}$, $i=0,1,\ldots,d$. Assume that $y_{1}, y_{2}, \ldots, y_{m-n}$ are all unit orthogonal vectors such that $R_{G}y_k=0$. Then the signless Laplacian spectrum of $\mathcal{Q}$-graph of $G$  consists precisely of the following:
\begin{enumerate}[(i)]
	\item $2r-2$ is the signless Laplacian eigenvalue of $\mathcal{Q}$-graph, with multiplicity $m-n$ and the corresponding orthogonal eigenvectors are
	\[
	Y_k=\frac{1}{||{y_k}||}\left({\begin{array}{*{20}{c}}
	  0\\
	  {{y_k}}
	  \end{array}}\right),
	\]
	for $k=1,2,\ldots,m-n$.
	\item $ q_{i\pm}=\frac{3r+q_{i}-2 \pm \sqrt{(q_{i}+r-2)^{2}+4q_{i}}}{2}$ are the signless Laplacian eigenvalues of $\mathcal {Q}$-graph and the corresponding orthogonal eigenvectors are
	\[
	{X_{i \pm }^{j}} = \frac{1}{{\sqrt {{{({q_{i\pm}}+2-2r- {q_i})}^2+{q_i}}} }}\left( {\begin{array}{*{20}{c}}
		{({q_{i \pm }} + 2 - 2r - {q_i}){x_j^{(i)}}}\\
		{R_G^T{x_j^{(i)}}}
		\end{array}} \right)
	\]
	for $j=1,2,\ldots,l_i$ and $i=0,1,\ldots,d$.
\end{enumerate}

\begin{proof}
According to the definition of $\mathcal {Q}$-graph, then the signless Laplacian matrix of $\mathcal {Q}$-graph of $G$ is given by
\begin{equation*}
{Q_{\mathcal {Q}(G)}} = \left( {\begin{array}{*{20}{c}}
{r{I_n}}&{{R_G}}\\
{{R_G}^T}&{2r{I_m} + A(\ell (G))},
\end{array}} \right)
\end{equation*}
where $\ell(G)$ denotes the line graph of $G$. The proof this theorem is divided into three claims.

\paragraph{Claim 1.} $2r-2$ and $q_{i\pm}$ are the signless Laplacian eigenvalues of $\mathcal {Q}$-graph of $G$ with corresponding respective multiplicities $m-n$ and $l_{i}$ for $i=0, 1, 2, \ldots, d$.
\\\\
\textit{Proof of Claim $1$.} Recall that Claim 1 can be obtained easily from Theorem 5.6 in \cite{J.-P. Li}, here we give the detailed proof for the convenience of readers. Since $A_{\ell(G)}=R_{G}^{T}R_{G}-2I_{m}$ and $Q_{G}=R_{G}R_{G}^{T}$, then the signless Laplacian characteristic polynomial of $\mathcal {Q}(G)$
\begin{equation*}
\begin{split}
{P_{{Q_{\mathcal {Q}(G)}}}}(t) &=det \left( {\begin{array}{*{20}{c}}
{(t - r){I_n}}&{ - {R_G}}\\
{ - R_G^T}&{(t - 2r){I_m} - A(\ell (G))}
\end{array}} \right)\\
& = det\left( {\begin{array}{*{20}{c}}
{(t - r){I_n}}&{ - {R_G}}\\
{ - R_G^T}&{(t + 2 - 2r){I_m} - R_G^T{R_G}}
\end{array}} \right)\\
& =det \left( {\begin{array}{*{20}{c}}
{(t - r){I_n}}&{ - {R_G}}\\
{(r - t - 1)R_G^T}&{(t + 2 - 2r){I_m}}
\end{array}} \right)\\
&=det \left( {\begin{array}{*{20}{c}}
{(t - r){I_n} + \frac{{r - t - 1}}{{t - 2r + 2}}{R_G}R_G^T}&0\\
{(r - t - 1)R_G^T}&{(t + 2 - 2r){I_m}}
\end{array}} \right)\\
& = {(t + 2 - 2r)^m}\det ((t - r){I_n} + \frac{{r - t - 1}}{{t - 2r + 2}}{R_G}R_G^T)\\
& = {(t + 2 - 2r)^{m - n}}\det ((t + 2 - 2r)(t - r){I_n} + (r - t - 1){R_G}{R_G}^T)\\
&= {(t + 2 - 2r)^{m - n}}\det ((t + 2 - 2r)(t - r){I_n} + (r - t - 1){Q_G}).
\end{split}
\end{equation*}
Since the distinct signless Laplacian eigenvalues of $G$ are $q_{0}>q_{1}>\cdots>q_{d}$ with respective multiplicities $l_{0}, l_{1}, \ldots, l_{d}$. Then we have
\[
{P_{{Q_{\mathcal {Q}(G)}}}}(t) = {(t + 2 - 2r)^{m - n}} \prod \limits_{i = 0}^d {[(t + 2 - 2r)(t - r) + (r - t - 1){q_i}]^{{l_i}}},
\]
which implies that $2r-2$ and
\[
{q_{i \pm }} = \frac{{3r + {q_i} - 2 \pm \sqrt {{{(q_{i} + r - 2)}^2} + 4{q_i}} }}{2}
\]
are the signless Laplacian eigenvalues of $\mathcal {Q}(G)$ for $i=0,1,\ldots,d$. Since $G$ is a non-bipartite graph, then $q_d\neq 0$. This implies that $q_{i \pm }\neq2r-2$ for $i=0,1,\ldots,d$. Hence, $2r-2$ and $q_{i \pm }$ are the signless Laplacian eigenvalues of $\mathcal {Q}(G)$ with corresponding respective multiplicities $m-n$ and $l_{i}$ for $i=0, 1, 2, \ldots, d$.

\paragraph{Claim 2.} The $Y_k$ and $X_{i \pm }^j$, as described by the theorem, are the signless Laplacian eigenvectors corresponding the signless Laplacian eigenvalues $2r-2$ and $q_{i \pm }$, respectively.
\\\\
\textit{Proof of Claim $2$.} By a simple calculation, we have
\[
Q_{\mathcal {Q}(G)}Y_k=
\left( {\begin{array}{*{20}{c}}
{r{I_n}}&{{R_G}}\\
{{R_G}^T}&{2r{I_m} + A(\ell (G))}
\end{array}} \right)\frac{1}{||{y_k}||}\left( {\begin{array}{*{20}{c}}
0\\
y_k
\end{array}} \right) = (2r - 2)\frac{1}{||{y_k}||}\left( {\begin{array}{*{20}{c}}
0\\
y_k
\end{array}} \right)=(2r - 2)Y_k,
\]
which implies that, for every $y ={y_1}, {y_2}, \ldots, {y_{m - n}}$,
\[
Y_k = \frac{1}{||{y_k}||}\left( {\begin{array}{*{20}{c}}
0\\
y_k
\end{array}} \right)
\]
is the signless Laplacian eigenvector corresponding to the signless Laplacian eigenvalue $2r-2$. Since $\{x_1^{(i)}, x_2^{(i)}, \ldots, x_{{l_i}}^{(i)}\}$ is an orthogonal basis of the eigenspace $V_{q_{i}}$ corresponding to the signless Laplacian eigenvalue $q_{i}$. Then $Q_{G}x_{j}^{(i)}=q_{i}x_{j}^{(i)}$ for $j=1, 2, \ldots, l_{i}$ and $i=0, 1, 2, \ldots, d$. It is easy to see that
\[
{q_{i \pm }}= r + \frac{{{q_i}}}{{{q_{i \pm }} + 2 - 2r - {q_i}}}.
\]
Now we obtain
\begin{equation*}
\begin{split}
{Q_{\mathcal {Q}(G)}}X_{i \pm }^j&=\frac{1}{{\sqrt {{{({q_{i \pm }} + 2 - 2r - {q_i})}^2} + {q_i}} }}\left( {\begin{array}{*{20}{c}}
{r{I_n}}&{{R_G}}\\
{{R_G}^T}&{2r{I_m} + A(\ell (G))}
\end{array}} \right)\left( {\begin{array}{*{20}{c}}
{({q_{i \pm }} + 2 - 2r - {q_i})x_j^{(i)}}\\
{R_G^Tx_j^{(i)}}
\end{array}} \right)\\
&= \frac{1}{{\sqrt {{{({q_{i \pm }} + 2 - 2r - {q_i})}^2} + {q_i}} }}\left( {\begin{array}{*{20}{c}}
{r{I_n}}&{{R_G}}\\
{{R_G}^T}&{(2r - 2){I_m} + R_G^T{R_G}}
\end{array}} \right)\left( {\begin{array}{*{20}{c}}
{({q_{i \pm }} + 2 - 2r - {q_i})x_j^{(i)}}\\
{R_G^Tx_j^{(i)}}
\end{array}} \right)\\
&= \frac{1}{{\sqrt {{{({q_{i \pm }} + 2 - 2r - {q_i})}^2} + {q_i}} }}\left( {\begin{array}{*{20}{c}}
{[r({q_{i \pm }} + 2 - 2r - {q_i}) + {q_i}]x_j^{(i)}}\\
{({q_{i \pm }} + 2 - 2r - {q_i} + 2r - 2 + {q_i})R_G^Tx_j^{(i)}}
\end{array}} \right)\\
&= \frac{1}{{\sqrt {{{({q_{i \pm }} + 2 - 2r - {q_i})}^2} + {q_i}} }}{q_{i \pm }}\left( {\begin{array}{*{20}{c}}
{({q_{i \pm }} + 2 - 2r - {q_i})x_j^{(i)}}\\
{R_G^Tx_j^{(i)}}
\end{array}} \right)\\
&= {q_{i \pm }}X_{i \pm }^j.
\end{split}
\end{equation*}
Hence, $X_{i\pm}^{j}$ are the signless Laplacian eigenvector of $\mathcal {Q}(G)$ corresponding to the signless Laplacian eigenvalues $q_{i\pm}$.

\paragraph{Claim 3.} All $X_{i\pm}^{j}$'s and $Y_k$'s are orthogonal signless Laplacian eigenvectors of $\mathcal {Q}(G)$.
\\\\
\textit{Proof of Claim $3$.} According to Claim $2$, we have $(X_{i\pm}^{j})^TX_{i\pm}^{k}=0$ and $Y_j^{T}Y_k=0$ for $j\neq k$. It is easy to know that $(X_{i\pm}^{j})^TY_{k}=0$ for $j=1, 2, \ldots, l_{i}$, $i=0, 1, 2, \ldots, d$ and $k=1,2,\ldots, m-n$.

In what follows, we will only prove that $X_{i+ }^jX_{i- }^j=0$. For the sake of convenience, set
$$\alpha=\frac{1}{{\sqrt {{{({q_{i + }} + 2 - 2r - {q_i})}^2+{q_i}}} }}\frac{1}{{\sqrt {{{({q_{i- }} + 2 - 2r - {q_i})}^2+{q_i}}} }}.$$
Since $({q_{i + }} + 2 - 2r - {q_i})({q_{i - }} + 2 - 2r - {q_i}) =  - {q_i}$, then
\begin{equation*}
\begin{split}
{(X_{i + }^j)^T}X_{i - }^j &= \alpha{\left( {\begin{array}{*{20}{c}}
{({q_{i _{+} }} + 2 - 2r - {q_i})x_j^{(i)}}\\
{R_G^Tx_j^{(i)}}
\end{array}} \right)^T}\left( {\begin{array}{*{20}{c}}
{({q_{i _{-} }} + 2 - 2r - {q_i})x_j^{(i)}}\\
{R_G^Tx_j^{(i)}}
\end{array}} \right)\\
&= \alpha(({q_{i + }} + 2 - 2r - {q_i})({q_{i - }} + 2 - 2r - {q_i}){(x_j^{(i)})^T}x_j^{(i)} + {(x_j^{(i)})^T}{R_G}R_G^Tx_j^{(i)})\\
&=  \alpha(- {q_i}{(x_j^{(i)})^T}x_j^{(i)} + {q_i}{(x_j^{(i)})^T}x_j^{(i)})\\
&= 0.
\end{split}
\end{equation*}
Therefore, $X_{i+}^j$ and $X_{i - }^j$ are orthogonal eigenvectors for any $j=1, 2, \ldots, l_{i}$ and $i=0, 1, 2, \ldots, d$.
Above all, all $X_{i\pm}^j$'s and $Y_k$'s are orthogonal eigenvectors for  $j=1, 2, \ldots, l_{i}$, $i=0, 1, \ldots, d$ and $k=1,2,\ldots, m-n$.
\end{proof}

\paragraph{Theorem 3.2.}
Let $G$ be an $r$-regular $(r\geq2)$ bipartite connected graph of order $n$ with $m$ edges, and $V_{1}\cup V_{2}$ be the bipartition of the vertex set of $G$. Also let $ 2r=q_{0} > q_{1}> \cdots > q_{d}=0$ are different signless Laplacian eigenvalues of $G$ with the corresponding multiplicities $l_{0}, l_{1}, \ldots l_{d}$, and $x_{1}^{(i)},x_{2}^{(i)}, \ldots, x_{l_{i}}^{(i)}$ be the unit orthogonal eigenvectors corresponding to eigenvalue $q_{i}$ for $i=0,1,\ldots,d$. Suppose that $y_{1}, y_{2}, \ldots, y_{m-n+1}$ are all unit orthogonal vectors such that $R_{G}y_k=0$ for $k=1,2,\ldots,m-n+1$. Then $2r-2$, $r$ and
$$ q_{i\pm}=\frac{3r+q_{i}-2 \pm \sqrt{(q+r-2)^{2}+4q_{i}}}{2}$$
are the signless Laplacian eigenvalues of $\mathcal {Q}$-graph and the corresponding signless Laplacian eigenvectors are, for $k=1,2,\ldots, m-n+1$,
\[
Y_k = \frac{1}{||{y_k}||}\left( {\begin{array}{*{20}{c}}
	0\\
	{{y_k}}
	\end{array}} \right),\;\;
\frac{1}{{\sqrt n }}\left( {\begin{array}{*{20}{c}}
	{{j_{|{V_1}|}}}\\
	{ - {j_{|{V_2}|}}}\\
	{{0_m}}
	\end{array}} \right)
\]
and
\[
{X_{i \pm }^{j}} = \frac{1}{{\sqrt {{{({q_{i \pm }} + 2 - 2r - {q_i})}^2+{q_i}}} }}\left( {\begin{array}{*{20}{c}}
{({q_{i \pm }} + 2 - 2r - {q_i}){x_j^{(i)}}}\\
{R_G^T{x_j^{(i)}}}
\end{array}} \right)
\]
for $j=1, 2, \ldots, l_{i}$, $i=0, 1, 2, \ldots, d-1$.

\begin{proof}
Since $G$ is bipartite, then the smallest signless Laplacian eigenvalue $q_{d}=0$. In the light of
\[
{q_{i \pm }} = \frac{{3r + {q_i} - 2 \pm \sqrt {{{(q_{i} + r - 2)}^2} + 4{q_i}} }}{2},
\]
we have ${q_{d + }}= 2r - 2$ and ${q_{d - }} = r$.
It is easy to verify that
\[
{Q_{\mathcal {Q}(G)}}\frac{1}{{\sqrt n }}\left( {\begin{array}{*{20}{c}}
{{j_{|{V_1}|}}}\\
{ - {j_{|{V_2}|}}}\\
{{0_m}}
\end{array}} \right)\\
 = \left( {\begin{array}{*{20}{c}}
{r{I_n}}&{{R_G}}\\
{{R_G}^T}&{(2r - 2){I_m} + R_G^T{R_G}}
\end{array}} \right)\frac{1}{{\sqrt n }}\left( {\begin{array}{*{20}{c}}
{{j_{|{V_1}|}}}\\
{ - {j_{|{V_2}|}}}\\
{{0_m}}
\end{array}} \right)\\
 = r\frac{1}{{\sqrt n }}\left( {\begin{array}{*{20}{c}}
{{j_{|{V_1}|}}}\\
{ - {j_{|{V_2}|}}}\\
{{0_m}}
\end{array}} \right).
\]
Thus
\[\frac{1}{{\sqrt n }}\left( {\begin{array}{*{20}{c}}
{{j_{|{V_1}|}}}\\
{ - {j_{|{V_2}|}}}\\
{{0_m}}
\end{array}} \right)\]% MathType!End!2!1!
is the signless Laplacian eigenvector corresponding to the signless Laplacian eigenvalue $r$ of $\mathcal {Q}(G)$. The rest of the proof is similar exactly to that of Theorem 3.1, omitted.
\end{proof}

According to Theorems 3.1 and 3.2, we obtain immediately the signless Laplacian eigenprojectors and spectral decomposition of the signless Laplacian matrix of the $\mathcal {Q}$-graph of $G$.

\paragraph{Theorem 3.3.} Assume that $G$ is an $r$-regular $(r\geq2)$ connected graph of order $n$ with $m$ edges. Then
\begin{enumerate}[(a)]
\item If $G$ is non-bipartite, then $F_{q_{i\pm}}$ and $F_{2r-2}$ are the  eigenprojectors corresponding to the respective signless Laplacian eigenvalues $q_{i\pm}$ and $2r-2$ of $\mathcal {Q}(G)$, where
\begin{equation}\small\label{9}
{F_{{q_{i \pm }}}} = \frac{1}{{{{({q_{i \pm }} + 2 - 2r - {q_i})}^2} + {q_i}}}\left( {\begin{array}{*{20}{c}}
{{{({q_{i \pm }} + 2 - 2r - {q_i})}^2}{f_{{q_i}}}}&{({q_{i \pm }} + 2 - 2r - {q_i}){f_{{q_i}}}{R_G}}\\
{({q_{i \pm }} + 2 - 2r - {q_i}){{({f_{{q_i}}}{R_G})}^T}}&{R_G^T{f_{{q_i}}}{R_G}}
\end{array}} \right)
\end{equation}
and
\begin{equation}\label{10}
{F_{2r - 2}} = \sum\limits_{k = 1}^{m - n} {\frac{1}{{||{y_k}||^{2}}}} \left( {\begin{array}{*{20}{c}}
0&0\\
0&{{y_k}y_k^T}
\end{array}} \right).
\end{equation}
Thus, we get the spectral decomposition of $Q_{\mathcal {Q}(G)}$ as follows:
\begin{equation}
{Q_{\mathcal {Q}(G)}} = \sum\limits_{i = 0}^{d} {\sum\limits_ \pm  {{q_{i \pm }}{F_{{q_{i \pm }}}} + (2r - 2)} } {F_{2r - 2}}.
\end{equation}
\item If $G$ is bipartite, then $F_{2r-2}$, $F_{r}$ and $F_{q_{i\pm}}$ are  the eigenprojectors corresponding to the respective signless Laplacian eigenvalues $2r-2$, $r$ and $q_{i\pm}$ of $\mathcal {Q}(G)$, where
\begin{equation}
{F_{2r - 2}} = \sum\limits_{k = 1}^{m - n+1} {\frac{1}{{||{y_k}||^{2}}}} \left( {\begin{array}{*{20}{c}}
0&0\\
0&{{y_k}y_k^T}
\end{array}} \right),
\end{equation}
\begin{equation}
{F_r} = \frac{1}{n}\left( {\begin{array}{*{20}{c}}
{{J_{|{V_1}|}}}&{ - {J_{|{V_1}| \times |{V_2}|}}}&0\\
{ - {J_{|{V_2}| \times |{V_1}|}}}&{{J_{|{V_2}|}}}&0\\
0&0&0
\end{array}} \right) = \left( {\begin{array}{*{20}{c}}
{{f_0}}&0\\
0&0
\end{array}} \right),
\end{equation}
and for $i\neq d$
\begin{equation}\small
{F_{{q_{i \pm }}}} = \frac{1}{{{{({q_{i \pm }} + 2 - 2r - {q_i})}^2} + {q_i}}}\left( {\begin{array}{*{20}{c}}
	{{{({q_{i \pm }} + 2 - 2r - {q_i})}^2}{f_{{q_i}}}}&{({q_{i \pm }} + 2 - 2r - {q_i}){f_{{q_i}}}{R_G}}\\
	{({q_{i \pm }} + 2 - 2r - {q_i}){{({f_{{q_i}}}{R_G})}^T}}&{R_G^T{f_{{q_i}}}{R_G}}
	\end{array}} \right).
\end{equation}
Thus, we get the spectral decomposition of $Q_{\mathcal {Q}(G)}$ as follows:
\begin{equation}
{Q_{\mathcal {Q}(G)}} =\sum\limits_{i = 0}^{d - 1} {\sum\limits_ \pm  {{q_{i \pm }}{F_{{q_{i \pm }}}} + (2r - 2)} } {F_{2r - 2}} + r{F_r}.
\end{equation}
\end{enumerate}

In accordance with Theorems 3.1, 3.2 and 3.3, we obtain easily the following proposition, which will be applied to analyse the signless Laplacian state transfer of $\mathcal {Q}(G)$.

\paragraph{Proposition 3.4.} Let $G$ be an $r$-regular $(r\geq2)$ connected graph of order $n$ with $m$ edges. Then, for two vertices $u$ and $v$ of $G$, we have
\begin{enumerate}[(i)]
\item If $G$ is non-bipartite, then
\begin{equation}\small
{(e_v^{m + n})^T}\text{exp}( - it{Q_{\mathcal {Q}(G)}})e_u^{m + n}
 = {e^{ - it\frac{{3r - 2}}{2}}}\sum\limits_{i = 0}^{d} {{e^{ - it\frac{{{q_i}}}{2}}}} e_v^T{f_{{q_i}}}{e_u}(\cos\frac{{{\Delta _{{q_i}}}t}}{2} + \frac{{{q_i} + r - 2}}{{{\Delta _{{q_i}}}}}\sin\frac{{{\Delta _{{q_i}}}t}}{2}),
\end{equation}
where $\Delta_{q_i}=\sqrt{(q_i+r-2)^2+4q_i}$ for $i=0,1,\ldots,d$.
\item If $G$ is bipartite, then
\begin{equation}\small
\begin{split}
(e_v^{m + n})&^T\text{exp}( - it{Q_{\mathcal {Q}(G)}})e_u^{m + n}\\
 &= {e^{ - it\frac{{3r - 2}}{2}}}\sum\limits_{i = 0}^{d - 1} {{e^{ - it\frac{{{q_i}}}{2}}}} e_v^T{f_{{q_i}}}{e_u}(\cos\frac{{{\Delta _{{q_i}}}t}}{2} + \frac{{{q_i} + r - 2}}{{{\Delta _{{q_i}}}}}\sin\frac{{{\Delta _{{q_i}}}t}}{2}) + {e^{ - itr}}e_v^T{f_0}{e_u}.
\end{split}
\end{equation}
\end{enumerate}

\begin{proof} (i) Since $G$ is a non-bipartite graph. Then the (a) of Theorem 3.3 and (\ref{7}) imply that
\[
\text{exp}( - it{Q_{\mathcal {Q}(G)}}) = \sum\limits_{i = 0}^{d} {\sum\limits_ \pm  {{e^{ - it{q_{i \pm }}}}{F_{{q_{i \pm }}}} + {e^{ - it(2r - 2)}}} } {F_{2r - 2}}.
\]
By a simple calculation, we get
$$({q_{i + }} + 2 - 2r - {q_i})({q_{i - }} + 2 - 2r - {q_i}) =  - {q_i},$$
$${({q_{i + }} + 2 - 2r - {q_i})^2} + {({q_{i - }} + 2 - 2r - {q_i})^2} = {\Delta _{{q_i}}}^2 - 2{q_i}$$
and
$${({q_{i - }} + 2 - 2r - {q_i})^2} - {({q_{i + }} + 2 - 2r - {q_i})^2} = ({q_i} + r - 2){\Delta _{{q_i}}}.$$
Thus, from the formulas (\ref{9}) and (\ref{10}), we obtain
\begin{equation*}
\begin{split}
{(e_v^{m + n})^T}\text{exp}( - it{Q_{\mathcal {Q}(G)}})e_u^{m + n}
 &= {(e_v^{m + n})^T}(\sum\limits_{i = 0}^{d} {\sum\limits_ \pm  {{e^{ - it{q_{i \pm }}}}{F_{{q_{i \pm }}}} + {e^{ - it(2r - 2)}}} } {F_{2r - 2}})e_u^{m + n}\\
 &= \sum\limits_{i = 0}^{d} {\sum\limits_ \pm  {{e^{ - it{q_{i \pm }}}}{{(e_v^{m + n})}^T}{F_{{q_{i \pm }}}}} } e_u^{m + n}\\
 &= {e^{ - it\frac{{3r - 2}}{2}}}\sum\limits_{i = 0}^{d} {{e^{ - it\frac{{{q_i}}}{2}}}} e_v^T{f_{{q_i}}}{e_u}(\cos\frac{{{\Delta _{{q_i}}}t}}{2} + \frac{{{q_i} + r - 2}}{{{\Delta _{{q_i}}}}}\sin\frac{{{\Delta _{{q_i}}}t}}{2})
\end{split}
\end{equation*}
for two vertices $u$ and $v$ of $G$.

(ii) Assume that $G$ is a bipartite graph. It follows from the (b) of Theorem 3.3 and (\ref{7}) that
\[
\text{exp}( - it{Q_{\mathcal {Q}(G)}}) = \sum\limits_{i = 0}^{d} {\sum\limits_ \pm  {{e^{ - it{q_{i \pm }}}}{F_{{q_{i \pm }}}} + {e^{ - it(2r - 2)}}} } {F_{2r - 2}} + {e^{ - itr}}{F_r}.
\]
The rest proof is exactly similar to that of the (i), we omit the process of the proof.
\end{proof}

\section{Signless Laplacian perfect state transfer of $\mathcal {Q}$-graph}

In this section, we mainly study the existence of signless Laplacian perfect state transfer of $\mathcal {Q}$-graph. First we give a key lamma, which is similar to the Lamma $4.1$ in \cite{Yipeng Li}.

\paragraph{Lemma 4.1.} Let $G$ be an $r$-regular connected graph of order $n$ with $m$ edges and $r\geq2$. Also let $2r=q_{0}>q_{1}>\cdots>q_{d}$ be all distinct signless Laplacian eigenvalues of $G$, with the signless Laplacian eigenprojectors $f_{q_{0}}, f_{q_{1}},\ldots,f_{q_{d}}$, respectively.
\begin{enumerate}[(i)]
\item If $G$ is non-bipartite, then there is some index $i_0\in\{1,2,\ldots,d\}$ such that $(f_{q_{i_0}}R_{G})e_{k}^{m}\neq0$ for any $k\in\{1, 2,\ldots, m\}$.

\item If $G$ is bipartite, then there is some index $i_0\in\{1,2,\ldots, d-1\}$ such that $(f_{q_{i_0}}R_{G})e_{k}^{m}\neq0$ for any $k\in\{1, 2,\ldots, m\}$.
\end{enumerate}

\begin{proof} (i) Since $G$ is an $r$-regular non-bipartite graph, then the signless Laplacian eigenprojector $f_{q_{0}}=f_{2r}=\frac{1}{n}J_{n}$. In the light of $\sum\nolimits_{i = 0}^d{{f_{{q_i}}}}=I_{n}$, we have
$\sum\nolimits_{i = 1}^d {{f_{{q_i}}}}  = {I_n} - {f_{2r}} = {I_n} - \frac{1}{n}{J_n}$.
Thus,
\begin{equation*}
\sum\limits_{i = 1}^d {{f_{{q_i}}}}R_{G}e_{k}^{m}=({I_n} - \frac{1}{n}{J_n})R_{G}e_{k}^{m}\neq0,
\end{equation*}
which implies that, for any $k\in\{1, 2,\ldots, m\}$, there exists some index $i_0\in\{1,2,\ldots,d\}$ such that $(f_{q_{i_0}}R_{G})e_{k}^{m}\neq0$.

(ii) In this case, assume that $V_{1}\cup V_{2}$ is the bipartition of the vertex set of $G$. Since $G$ is a regular graph, then $|V_{1}|=| V_{2}|=\frac{n}{2}$. According to $\sum\nolimits_{i = 0}^d {{f_{{q_i}}}}  = I_{n}$, $f_{q_{0}}=f_{2r}=\frac{1}{n}J_{n}$ and
\[
{f_{{q_d}}} = {f_0} = \left( {\begin{array}{*{20}{c}}
{{J_{\frac{n}{2}}}}&{ - {J_{\frac{n}{2}}}}\\
{ - {J_{\frac{n}{2}}}}&{{J_{\frac{n}{2}}}}
\end{array}} \right),
\]
we obtain
\begin{equation*}
\sum\limits_{i = 1}^{d-1}{{f_{{q_i}}}}  = {I_n} - {f_{2r}}-{f_{0}} = {I_n} - \frac{1}{n}{J_n}-\left( {\begin{array}{*{20}{c}}
{{J_{\frac{n}{2}}}}&{ - {J_{\frac{n}{2}}}}\\
{ - {J_{\frac{n}{2}}}}&{{J_{\frac{n}{2}}}}
\end{array}} \right).
\end{equation*}
Hence, for any $k\in\{1, 2,\ldots, m\}$,
\begin{equation*}
\sum\limits_{i = 1}^{d-1} {{f_{{q_i}}}}R_{G}e_{k}^{m}=\left({I_n} - \frac{1}{n}{J_n}-\left( {\begin{array}{*{20}{c}}
{{J_{\frac{n}{2}}}}&{ - {J_{\frac{n}{2}}}}\\
{ - {J_{\frac{n}{2}}}}&{{J_{\frac{n}{2}}}}
\end{array}} \right)\right)R_{G}e_{k}^{m}\neq0,
\end{equation*}
which implies that the required result follows.
\end{proof}

\paragraph{Theorem 4.2.} Let $G$ be an $r$-regular connected graph of order $n\geq2$ with $m$ edges.
If all the signless Laplacian eigenvalues of $G$ are integers, then there is no signless Laplacian perfect state transfer in $\mathcal {Q}(G)$.

\begin{proof} Here we only proof the case of non-bipartite graph $G$. For bipartite graph, the proof is similar to that of non-bipartite graph.

Let $V(G)\cup I(G)$ be the vertex set of $\mathcal {Q}(G)$, where $I(G)$ is the new vertex set in $\mathcal {Q}(G)$. Also let $z$ be a vertex of $\mathcal {Q}(G)$. Suppose towards the contradiction that $\mathcal {Q}(G)$ has signless Laplacian perfect state transfer between the vertex $z$ and another vertex. According to Theorem $2.1$, we only need to discuss it in the following two cases.\\

\emph{Case 1.} All the elements of $\text{supp}_{Q_{\mathcal {Q}(G)}}(z)$ are integers. In this case, we consider the following two subcases: the vertex $z\in V(G)$ and $z\in I(G)$.\\

\emph{Subcase 1.1.} The vertex $z\in V(G)$. Since $G$ is a connected graph, then there is a positive eigenvalue $q\in \text{supp}_{Q_{G}}(z)$. Thus $f_{q}e_{z}\neq 0$. It follows from Lemma 4.1 and the (a) of Theorem 3.3 that $F_{q\pm}e_{z}^{m+n}\neq 0$. This implies that $q_{\pm}\in\text{supp}_{Q_{\mathcal {Q}(G)}}(z)$. Recall that all of $q$, $q_{+}$ and $q_{-}$ are integers. Since
\[
{q_ + } = \frac{{3r + q - 2 + \sqrt {{{(r + q - 2)}^2} + 4q} }}{2},\;\;
{q_ - } = \frac{{3r + q - 2 - \sqrt {{{(r + q - 2)}^2} + 4q} }}{2}.
\]
Then $q_{+}+q_{-}=3r+q-2$, $q_{+}-q_{-}=\sqrt{(r+q-2)^{2}+4q}$. It follows that $\sqrt{(r+q-2)^{2}+4q}$ is integer, equivalently, $(r+q-2)^{2}+4q$ is a perfect square. Clearly, $4q$ is even and
\[
{(q + r - 2)^2} < {(q + r - 2)^2} + 4q < {(q + r - 2)^2} + 4(q + r) < {(q + r + 2)^2}.
\]
If $(q + r - 2)^2 + 4q \neq (r+q)^{2}$, then at most one of $q_{+}$ and $q_{-}$ is integer, which contradicts that all the elements of $\text{supp}_{Q_{\mathcal {Q}(G)}}(z)$ are integers. Otherwise, $(q + r - 2)^2 + 4q = (r+q)^{2}$. Thus, we have $r=1$, which implies that $G=K_2$ and
$\mathcal {Q}(G)=P_3$. It is well known\cite{Coutinho2015} that $P_3$ has no signless Laplacian
state transfer, which contradicts to our hypothesis.\\

\emph{Subcase 1.2.} The vertex $z\in I(G)$. According to Lemma 4.1, there exists some signless Laplacian eigenvalue $q$ such that $(f_{q}R_{G})e_{z}^{m}\neq0$. Then, according to the (a) of Theorem 3.3, we obtain that
\[
{F_{{q_{\pm} }}}e_z^{m + n} = \frac{1}{{{{({q_{\pm} } + 2 - 2r - q)}^2} + q}}\left( {\begin{array}{*{20}{c}}
{({q_{\pm} } + 2 - 2r - q)({f_q}{R_G})e_z^m}\\
{(R_G^T{f_q}{R_G})e_z^m}
\end{array}} \right)\neq0.
\]
Thus, $q_\pm\in \text{supp}_{Q_{\mathcal {Q}(G)}}(z)$ and $q_{\pm}$ are integer. Similar to the proof of Subcase $1.1$ above, we get that not all the elements of $\text{supp}_{Q_{\mathcal {Q}(G)}}(z)$ are integers, a contradiction.\\

\emph{Case 2.} All the elements of $\text{supp}_{Q_{\mathcal {Q}(G)}}(z)$ are quadratic integers. At this time, it is divided into two subcases below.\\

\emph{Subcase 2.1.} The vertex $z \in V(G)$. On the basis of Subcase 1.1, we have $ q_\pm\in \text{supp}_{Q_{\mathcal {Q}(G)}}(z)$ whenever $q\in \text{supp}_{Q_{G}}(z)$. Since all the elements of $\text{supp}_{Q_{\mathcal {Q}(G)}}(z)$ are quadratic integers and $(q + r - 2)^2 + 4q $ can not be a perfect square for integer $q$. Then, there exists a square-free integer $\Delta>1$, along with integers $a$ (we allow that $a=0$), $b_{+}$ and $b_{-}$, such that
\[
{q_\pm } = \frac{{a + {b_ \pm }\sqrt \Delta  }}{2}.
\]
Noting that $(q_{+} + 2 - 2r - q)({q_{-}} + 2 - 2r - q) = - q$. Thus,
\begin{equation*}
\frac{1}{4}{b_ + }{b_ - }\Delta  + \frac{1}{4}{(a + 4 - 4r - 2q)^2} + \frac{1}{4}\sqrt \Delta  ({b_ + } + {b_ - })(a + 4 - 4r - 2q) =  - q.
\end{equation*}
Since $\sqrt{\Delta}$ is irrational, then either $b_ +  + b_ - =0$ or $a+4-4r-2q=0$. If $b_ +  + b_ - =0$, then $q_{+}+q_{-}=a=3r+q-2$. At this time, $|\text{supp}_{Q_{G}}(z)|=1$, which contradicts that $G$ is a connected graph with $n\geq2$ vertices. Otherwise, $a+4-4r-2q=0$, then $a=4r+2q-4$, which also implies that $|\text{supp}_{Q_{G}}(z)|=1$, a contradiction. Hence, if all the elements of $\text{supp}_{Q_{\mathcal {Q}(G)}}(z)$ are quadratic integer, then neither $b_ +  + b_ - =0$ nor $a+4-4r-2q=0$.\\

\emph{Subcase 2.2.} The vertex $z\in I(G)$. According to Lamma $4.1$, there exists some signless Laplacian eigenvalue $q$ such that $(f_{q}R_{G})e_{z}^{m}\neq0$. Then, from Theorem $3.3$, we obtain that
\[
{F_{{q_\pm }}}e_z^{m + n} = \frac{1}{{{{({q_ \pm } + 2 - 2r - q)}^2} + q}}\left( {\begin{array}{*{20}{c}}
{({q_ \pm } + 2 - 2r - q)({f_q}{R_G})e_z^m}\\
{(R_G^T{f_q}{R_G})e_z^m}
\end{array}} \right)\neq0.
\]
So, $q_\pm \in \text{supp}_{Q_{\mathcal {Q}(G)}}(z)$ and $q_{\pm}$ are quadratic integers. Similar to the proof of Subcase 2.1 above, not all the elements of $\text{supp}_{Q_{\mathcal {Q}(G)}}(z)$ are quadratic integers.

Neither two cases above is likely to happen. Hence, it follows from Theorem $2.1$ that there is no signless Laplacian perfect state transfer in $\mathcal {Q}(G)$.
\end{proof}

Remark that the path $P_{2}$ is signless Laplacian integral as the signless Laplacian spectrum of $P_{2}$ is $\{0,2\}$. Theorem 4.2 shows that its $\mathcal {Q}$-graph $P_{3}$ has no signless laplacian perfect state transfer. This result has been proved by Coutinho and Liu in \cite{Coutinho2015}. In addition, in accordance with the proof of Theorem 4.2, we obtain easily the following corollary.

\paragraph{Corollary 4.3.} Let $G$ be an $r$-regular connected graph of order $n\geq2$ with $m$ edges. Assume that $u$ is any vertex in $G$.
If all the elements of $\text{supp}_{Q_G}(u)$ are integers, then $\mathcal {Q}(G)$ has no signless Laplacian perfect state transfer between $u$ and any other vertex of $\mathcal {Q}(G)$.

\paragraph{Example 1.} Assume that $G$ is one of the following graph:
\begin{enumerate}[(i)]
\item $d$-dimensional hypercubes $Q_{d}$ for $d\geq2$;
\item Cocktail party graphs $\overline{mK_{2}}$ for $m\geq 2$;
\item Halved $2d$-dimensional hypercubes $\frac{1}{2}Q_{2d}$ for $d\geq2$, where the vertex set $V(\frac{1}{2}Q_{2d})$ consists of elements of $\mathbb{Z}_{2}^{2d}$ of even Hamming weight and two vertices are adjacent if and only if their Hamming distance is exactly two.
\end{enumerate}
Then $\mathcal {Q}(G)$ has no signless Laplacian perfect state transfer.

\begin{proof}
It has been showed \cite{Tian2021} that the signless Laplacian spectra of these graphs are given by
 \begin{enumerate}[(i)]
\item $\text{Sp}(Q_{d})=\{2d-2l:0\leq l\leq d\}$;
\item $\text{Sp}(\overline{mK_{2}})=\{4m-4, 2m-2, 2m-4\}$;
\item $\text{Sp}(\frac{1}{2}Q_{2d})=\{2
\left( {\begin{array}{*{20}{c}}
{2d}\\
2
\end{array}} \right)-2l(2d-l): 0\leq l\leq d
\}$.
\end{enumerate}
Obviously, these graphs are regular connected graphs and all the signless Laplacian eigenvalues of these graphs are integers. Hence, if $G$ is each one of these graphs, then $\mathcal {Q}(G)$ has no signless Laplacian perfect state transfer by Theorem 4.2.
\end{proof}

\section{Signless Laplacian pretty good state transfer of $\mathcal{Q}$-graph}

Recall that the graph possessing Signless Laplacian perfect state transfer is rarely. In this section, we will discuss the existence of Signless Laplacian pretty good state transfer of $\mathcal {Q}$-graph and give some examples existing signless Laplacian pretty good state transfer.

\paragraph{Theorem 5.1} Let $G$ be an $r$-regular connected graph of order $n$ with $m$ edges and $r\geq2$. Also let $ q_{0}>q_{1}>\cdots>q_{d}$ be all different signless Laplacian eigenvalues of $G$ and $g$ is as described in Theorem 2.1. Assume that $G$ has signless Laplacian perfect state transfer between vertices $u$ and $v$. If $r$ is divisible by $g$, then $\mathcal {Q}(G)$ has signless Laplacian pretty good state transfer between vertices $u$ and $v$.

\begin{proof} The proof of this theorem is divided into the following two cases. In what follows, for the sake of convenience, set $S=\text{supp}_{Q_{G}}(u)$ and $\Delta_{q_i}=\sqrt{(q_i+r-2)^2+4q_i}$ for $i=0,1,\ldots,d$.\\

\emph{Case 1.} $G$ is a connected non-bipartite graph. In this case, from the (i) of Proposition 3.4, we have
\begin{equation*}
\begin{split}
{(e_v^{m + n})^T}\exp( - it{Q_{\mathcal {Q}(G)}})e_u^{m + n}
 = {e^{ - it\frac{{3r - 2}}{2}}}\sum\limits_{i = 0}^{ d} {{e^{ - it\frac{{{q_i}}}{2}}}} e_v^T{f_{{q_i}}}{e_u}\left(\cos\frac{{{\Delta _{{q_i}}}t}}{2} + \frac{{{q_i} + r - 2}}{{{\Delta _{{q_i}}}}}\sin\frac{{{\Delta _{{q_i}}}t}}{2}\right).
\end{split}
\end{equation*}

\begin{equation*}
\begin{split}
{{{(e_v^{m + n})}^T}\exp( - it{Q_{\mathcal {Q}(G)}})e_u^{m + n}}
%\begin{array}{l}
 &= {e^{ - it\frac{{3r - 2}}{2}}}\sum\limits_{i = 0}^{d} {{e^{ - it\frac{{{q_i}}}{2}}}} e_v^T{f_{{q_i}}}{e_u}\left(\cos\frac{{{\Delta _{{q_i}}}t}}{2} + \frac{{{q_i} + r - 2}}{{{\Delta _{{q_i}}}}}\sin\frac{{{\Delta _{{q_i}}}t}}{2}\right)\\
%\end{array}\\{}\\
&={e^{ - it\frac{{3r - 2}}{2}}}\sum\limits_{{q_i} \in S} {{e^{ - it\frac{{{q_i}}}{2}}}} e_v^T{f_{{q_i}}}{e_u}\left(\cos\frac{{{\Delta _{{q_i}}}t}}{2} + \frac{{{q_i} + r - 2}}{{{\Delta _{{q_i}}}}}\sin\frac{{{\Delta _{{q_i}}}t}}{2}\right),
\end{split}
\end{equation*}
where the last equality holds because $e_v^T{f_{{q_i}}}{e_u}\neq0$ if and only if $q_{i} \in \text{supp}_{Q_{G}}(u)$.

Clearly, $|e^{-it\frac{3r-2}{2}}|=1$ for any $t$. In order to prove that $\mathcal {Q}(G)$ has signless Laplacian pretty good state transfer between the vertices $u$ and $v$, we only need to find a time $t_0$ such that
\begin{equation}
\left|\sum\limits_{{q_i} \in S}^{} {{e^{ - it\frac{{{q_i}}}{2}}}} e_v^T{f_{{q_i}}}{e_u}\left(\cos\frac{{{\Delta _{{q_i}}}t_0}}{2} + \frac{{{q_i} + r - 2}}{{{\Delta _{{q_i}}}}}\sin\frac{{{\Delta _{{q_i}}}t_0}}{2}\right)\right| \approx 1.
\end{equation}
Observe that $f_{2r}e_{u}\neq 0$ implies that $2r\in S$. Since $G$ has signless Laplacian perfect state transfer at $u$ and $v$. Then, it follows from Theorem 2.1 and Lemma 2.4 that all the elements of $S$ are integers. This implies that $\mathcal {Q}(G)$ has no signless Laplacian perfect state transfer at $u$ and $v$ by Corollary 4.3. Bear in mind that $G$ is connected, $|S|\geq 2$. Furthermore, assume that $2r\neq q_{j}\in S$ for some $j\in\{1,2,\ldots, d\}$. Then $q_{j}$ is integer and $f_{q_{j}}e_{u}\neq 0$, which implies that $q_{j\pm} \in \text{supp}_{Q_{\mathcal {Q}(G)}}(u)$. In what follows, we let $\Delta _{q_{j}}=a_{j}\sqrt{b_{j}}$ for each $q_{j}\in S$, where $a_{j}, b_{j}\in \mathbb{Z}^{+}$ and $b_{j}$ is the square-free part of $\Delta _{q_{j}}^{2}$. Since $(q_{j}+r-2)^{2}<(q_{j}+r-2)^{2}+4q_{j}<(q_{j}+r)^{2}$ for $r>1$. Then the disjoint union
\[
\cup \{ \sqrt {{b_j}} :{q_j} \in S,{q_j} > 0\}
\]
is linearly independent over $\mathbb{Q}$ by Lemma 2.3. By Theorem $2.2$, there exist integers $\alpha$ and $c_{j}$ for each $q_{j}\in S$, such that
\begin{equation}\label{19}
\alpha \sqrt {{b_j}}  - {c_j} \approx  - \frac{1}{{2g}}\sqrt {{b_j}}.
\end{equation}
If $b_{l}=b_{j}$ for two distinct eigenvalues $q_{l}, q_{j}\in S$, then $c_{l}=c_{j}$. Multipling by $4a_{j}$ in both sides of (\ref{19}), we obtain
\[
\Delta_{q_{j}}\approx \frac{{4{a_j}{c_j}}}{{4\alpha  + \frac{2}{g}}}.
\]
Take $t_0=(4\alpha+\frac{2}{g})\pi$, then
\[\cos \frac{{{\Delta _{{q_j}}}t_0 }}{2} \approx \cos \frac{{\frac{{4{a_j}{c_j}}}{{4\alpha  + \frac{2}{g}}}(4\alpha  + \frac{2}{g})\pi}}{2} = \cos 2{a_j}{c_j}\pi = 1.
\]
Hence, we have
\begin{equation*}
\begin{split}
\left|\sum\limits_{{q_j} \in S} {{e^{ - it_0 \frac{{{q_j}}}{2}}}} e_v^T{f_{{q_j}}}{e_u}(cos\frac{{{\Delta _{{q_j}}}t_0 }}{2} + \frac{{{q_j} + r - 2}}{{{\Delta _{{q_i}}}}}sin\frac{{{\Delta _{{q_j}}}t_0 }}{2})\right|
 &\approx \left|\sum\limits_{{q_j} \in S} {{e^{ - it_0 \frac{{{q_j}}}{2}}}} e_v^T{f_{{q_j}}}{e_u}\right|\\
& \approx \left|\sum\limits_{{q_j} \in S} {{e^{ - i\frac{\pi }{g}{q_j}}}} e_v^T{f_{{q_j}}}{e_u}\right|=1,
\end{split}
\end{equation*}
where the last equality holds by Theorem 2.1.\\

\emph{Case 2.} $G$ is a connected bipartite graph. According to the (ii) of Proposition 3.4 and Case 1, we only need to prove that
\begin{equation}
\left|\sum\limits_{{q_j} \in S\backslash 0} {{e^{ - it\frac{{{q_j}}}{2}}}} e_v^T{f_{{q_j}}}{e_u}(\cos\frac{{{\Delta _{{q_j}}}t}}{2} + \frac{{{q_j} + r - 2}}{{{\Delta _{{q_j}}}}}\sin\frac{{{\Delta _{{q_j}}}t}}{2}) + {e^{ - itr}}e_v^T{f_0}{e_u}\right| \approx 1
\end{equation}
for some $t$. Take $t_0=(4\alpha+\frac{2}{g})\pi$ again. Similar to the discussion in Case 1, we give
\begin{equation}\label{21}
\begin{split}
&\left|\sum \limits_{{q_j} \in S\backslash 0} {{e^{ - it_0 \frac{{{q_j}}}{2}}}} e_v^T{f_{{q_j}}}{e_u}(\cos\frac{{{\Delta _{{q_j}}}t_0 }}{2} + \frac{{{q_j} + r - 2}}{{{\Delta _{{q_j}}}}}\sin\frac{{{\Delta _{{q_j}}}t_0 }}{2})+ {e^{ - it_0 r}}e_v^T{f_0}{e_u}\right|\\
&\approx\left|\sum\limits_{{q_j} \in S\backslash 0} {{e^{ - it_0 \frac{{{q_j}}}{2}}}} e_v^T{f_{{q_j}}}{e_u} + {e^{ - i\frac{r}{g}2\pi }}e_v^T{f_0}{e_u}\right| .
\end{split}
\end{equation}
Since $r$ is divisible by $g$, then ${e^{ - i\frac{r}{g}2\pi }}=e^{-it_0\frac{0}{2}}$. It follows from (\ref{21}) that
\begin{equation*}
\begin{split}
 \left|\sum\limits_{{q_j} \in S\backslash 0} {{e^{ - it_0 \frac{{{q_j}}}{2}}}} e_v^T{f_{{q_j}}}{e_u} + {e^{ - i\frac{r}{g}2\pi }}e_v^T{f_0}{e_u}\right|
 &= \left|\sum\limits_{{q_j} \in S\backslash 0} {{e^{ - it_0 \frac{{{q_j}}}{2}}}} e_v^T{f_{{q_j}}}{e_u} + {e^{ - it_0 \frac{0}{2}}}e_v^T{f_0}{e_u}\right|\\
 &= \left|\sum\limits_{{q_j} \in S} {{e^{ - it_0 \frac{{{q_j}}}{2}}}} e_v^T{f_{{q_j}}}{e_u}\right|\\
 & = \left|\sum\limits_{{q_j} \in S} {{e^{ - i\frac{\pi }{g}{q_j}}}} e_v^T{f_{{q_j}}}{e_u}\right|
 = 1.
\end{split}
\end{equation*}
To sum up two case above, we obtain the required result.
\end{proof}

In what following, we present some families of $\mathcal {Q}$-graphs admitting signless Laplacian pretty good state transfer for distance regular graphs. Distance regular graphs have many beautiful combinatorial properties that play an important role in graph theory and combinatorial mathematics. In quantum walks of graphs, an interesting result is the following: \emph{the eigenvalue support of each vertex in a distance regular graph $G$ equals the set of all distinct eigenvalues of $G$} \cite{Coutinho2015}. Since, for a regular graph, its signless Laplacian matrix shares the same eigenprojectors with its adjacency matrix. Then this property also holds for the signless Laplacian matrices in regular graphs. On the other hand, for a regular graph $G$, there exists perfect state transfer between vertices $u$ and $v$ if and only if $G$ has signless Laplacian perfect state transfer between vertices $u$ and $v$
\cite{Alvir2016}. Thus, applying Theorem 5.1, we may give some examples having signless Laplacian pretty good state transfer.

\paragraph{Example 2.} For the $d$-dimensional hypercube $Q_{d}$, cocktail party graphs $\overline{mK_{2}}$ and halved $2d$-dimensional hypercubes $\frac{1}{2}Q_{2d}$, these graphs admit signless Laplacian perfect state transfer between antipodal vertices $u$ and $v$. However, Example 1 showed that the $\mathcal {Q}$-graphs of these graph have no signless Laplacian perfect state transfer. In contrast to these results, we find that
\begin{enumerate}[(i)]
  \item The $\mathcal {Q}$-graph of $Q_{d}$ admits signless Laplacian pretty good state transfer between $u$ and $v$ whenever $d$ is even. Indeed, it is easy to check that $g=\gcd\left(\left\{\dfrac{q_0-q_k}{\sqrt{\Delta}}\right\}_{k=0}^d\right)=2$ and $Q_{d}$ is $d$-regular. The required result follows by Theorem 5.1;
  \item Note that $\overline{mK_{2}}$ is $(2m-2)$-regular and $g=2$ (see Example 1). Hence, $\mathcal {Q}(\overline{mK_{2}})$ has signless Laplacian pretty good state transfer between $u$ and $v$;
  \item Similarly, if $\left({\begin{array}{*{20}{c}}
	  2d\\
	  2
\end{array}}\right)$ is divisible by 2, then $\mathcal {Q}(\frac{1}{2}Q_{2d})$ admits signless Laplacian pretty good state transfer between $u$ and $v$.
\end{enumerate}

\section{Concluding remarks}
 The paper focuses on perfect state transfer and pretty good state transfer in $\mathcal {Q}$-graphs of regular graphs relative to signless Laplacian. We have obtained the spectral decomposition of signless Laplacian matrix of $\mathcal {Q}$-graph of a regular graph. Furthermore, some sufficient conditions have been given about $\mathcal {Q}$-graphs of regular graphs admitting perfect state transfer and pretty good state transfer relative to signless Laplacian. We end this paper with the following problem: \emph{For a non-regular graph, determine whether or not the corresponding $\mathcal {Q}$-graph has perfect state transfer or pretty good state transfer relative to signless Laplacian.}\\
\\
\textbf{Acknowledgements} This work was in part supported by the National Natural Science Foundation of China (No. 11801521).


\begin{thebibliography}{99}
	
\bibitem{Ackelsberg2016} E. Ackelsberg, Z. Brehm, A. Chan, J. Mundinger, C. Tamon, Laplaican state transfer in coronas, Linear Algebra Appl. 506 (2016) 154-167.

\bibitem{Ackelsberg2017} E. Ackelsberg, Z. Brehm, A. Chan, J. Mundinger, C. Tamon, Quantum State Transfer in Coronas, Electron. J. Combin. 24(2) (2017) $\#$P2.24.

\bibitem{Alvir2016} R. Alvir, S. Dever, B. Lovitz, J. Myer, C. Tamon, Y. Xu, H. Zhan, Perfect state transfer in Laplacian quantum walk, J. Algebraic Combin. 43(4) (2016) 801-826.

\bibitem{Banchi2017} L. Banchi, G. Coutinho, C. Godsil, S. Severini, Pretty good state transfer in qubit chains-the Heisenberg Hamiltonian, J. Math. Phys. 58(3) (2017) 032202, 9 pp.

\bibitem{Bose} S. Bose, Quantum communication through an unmodulated spin chain, Phys. Rev. Lett. 91(20) (2003) 207901.

\bibitem{S.Bose} S. Bose, A. Casaccino, S. Mancini, S. Severini, Communication in $XYZ$ all-to-all quantum networks with a missing link, Int. J. Quantum Inf. 7(4) (2009) 713-723.

\bibitem{Christandl} M. Christandl, N. Datta, A. Ekert, A. Landahl, Perfect state transfer in quantum spin networks, Phys. Rev. Lett. 92(18) (2004) 187902.

\bibitem{Coutinho2014} G. Coutinho, Quantum State Transfer in Graphs, University of Waterloo, 2014, PhD thesis.

\bibitem{Coutinho2016} G. Coutinho, C. Godsil, Perfect state transfer in products and covers of graphs, Linear Multilinear Algebra 64(2) (2016) 235-246.

\bibitem{Coutinho2015} G. Coutinho, H. Liu, No Laplacian perfect state transfer in trees, SIAM J. Discrete Math. 29(4) (2015) 2179-2188.

\bibitem{Cvetkovic1995} D. Cvetkovi\'{c}, M. Doob, H. Sachs, Spectra of Graphs: Theory and Application, Johann Ambrosius Barth Verlag, Heidelberg-Leipzig, 1995.


\bibitem{Cvetkovic2010} D. Cvetkovi\'{c}, P. Rowlinson, S.K. Simi\'{c}, An Introduction to the thoery of Graph Spectra, first edition, Cambridge University Press, 2010.


\bibitem{Cvetkovic2009} D. Cvetkovi\'{c}, S.K. Simi\'{c}, Towards a spectral theory of graphs based on signless Laplacian I. Publ. Inst. Math. (Beograd) (N.S.) 85(99) (2009) 19-33.

\bibitem{Cvetkovic2011} D. Cvetkovi\'{c}, S.K. Simi\'{c}, Graph spectra in computer science, Linear Algebra Appl. 434 (2011) 1545-1562.

\bibitem{Farhi} E. Farhi, S. Gutmann, Quantum computation and decision trees, Phys. Rev. A 58 (1998) 915-928.

\bibitem{Ge Yang} Y. Ge, B. Greenberg, O. Perez, C. Tamon, Perfect state transfer, graph products and equitable partitions, Int. J. Quantum Inf. 9(3) (2011) 823-842.

\bibitem{Godsil2012} C. Godsil, State transfer on graphs, Discrete Math. 312(1) (2012) 129-147.

\bibitem{equation} C. Godsil, When can perfect state transfer occur? Electron. J. Linear Algebra 23 (2012) 877-890.

\bibitem{numbertheory} G.H. Hardy, E.M. Wright, An Introduction to the Theory of Numbers, fifth edition, Oxford University Press, 2000.

\bibitem{Y.Hou} Y. Hou, W.-C. Shiu, The spectrum of the edge corona of two graphs, Electron. J. Linear Algebra 20 (2010) 586-594.

\bibitem{potential} M. Kempton, G. Lippner, S.-T. Yau, Perfect state transfer on graphs with potential, Quantum Inf. Comput. 17(3-4) (2017) 303-327.

\bibitem{involution} M. Kempton, G. Lippner, S.-T. Yau, Pretty good quantum state transfer in symmetric spin networks via magnetic field, Quantum Inf. Process. 16(9) (2017) Paper No. 210, 23 pp.

\bibitem{Yipeng Li} Y. Li, X. Liu, S. Zhang, Laplacian state transfer in $Q$-graph, Appl. Math. Comput. 384 (2020) 125370, 11 pp.

\bibitem{J.-P. Li} J. Li, B. Zhou, Signless Laplacian characteristic polynomials of regular graph transformations, arXiv: 1303.5527vl, 2013.

\bibitem{Xiaogang Liu} X. Liu, Q. Wang, Laplacian state transfer in total graphs, Discrete Math. 344 (2021) 112139, 11 pp.

\bibitem{galoistheory} I. Richards, An application of Galois theory to elementary arithmetic, Adv. Math. 13(3) (1974) 268-273.

\bibitem{Tian2021} G.-X. Tian, P.-K. Yu, S.-Y. Cui, The signless Laplacian state transfer in coronas, Linear Multilinear Algebra 69(2) (2021) 278-295.

\bibitem{Dam2003} E.R. van Dam, W. Haemers, Which graphs are determined by their spectrum? Linear Algebra Appl. 373 (2003) 241-272.

\bibitem{Wang2021} J. Wang, X. Liu, Laplacian state transfer in edge complemented coronas, Discrete Appl. Math. 293 (2021) 1-14.

\end{thebibliography}
\end{document}